\documentclass[11pt]{article}
\usepackage{amsfonts}
\usepackage{amsmath}
\usepackage{amssymb}

\usepackage[ top=2cm, bottom=2cm]{geometry}
\newtheorem{thm}{Theorem}[section]
\newtheorem{rmk}[thm]{Remark}

\newtheorem{cor}[thm]{Corollary}

\newtheorem{prf}{Proof}

\def\eps{\varepsilon}
\def\P{\mathbb{P}}

\def\N{\mathbb{N}}

\def\C{\mathbb{C}}

\def\bo{\nl\phantom{a}\hfill $\Box$\nl}

\def\nl{\newline}

\newcommand{\beq}{\begin{equation}}\newcommand{\eeq}{\end{equation}}
\usepackage{authblk}

\title{\bf \Large  The average order of the M\"{o}bius function for Beurling primes}
\author[1,2]{Ammar Ali Neamah}
\author[1]{Titus W Hilberdink}

\affil[1]{\small Department of Mathematics, University of Reading, Whiteknights, PO Box 220, Reading RG6 6AX, UK; t.w.hilberdink@reading.ac.uk}
\affil[2]{ Faculty of Computer Science and Mathematics, University of Kufa, Najaf, Iraq; ammara.meamah@uokufa.edu.iq,  a.a.n.al-rammahi@pgr.reading.ac.uk}

\date{}

\begin{document}
 \maketitle

\vspace{-10mm}
\indent
\begin{abstract} 
\vspace{0.2in}
In this paper, we study the counting functions $\psi_\mathcal{P}(x)$, $N_\mathcal{P}(x)$ and $M_\mathcal{P}(x)$ of a generalized prime system $\mathcal{N}$. Here $M_\mathcal{P}(x)$ is the partial sum of the M\"{o}bius function over $\mathcal{N}$  not exceeding $x$. In particular, we study these when they are asymptotically well-behaved, in the sense that $\psi_{\cal{P}}(x) =  x+O({x^{ \alpha+\eps }})$, $N_{\cal{P}}(x) = \rho x+O({x^{ \beta+\eps }})$ and $ M_\mathcal{P}(x) = O(x^{\gamma+\eps})$, for some $\rho >0$ and $\alpha, \beta, \gamma<1$. We show that the two largest of $\alpha,\beta,\gamma$ must be equal and at least $\frac{1}{2}$.

\

\vspace{0.2in}\end{abstract}

\noindent
{\bf 1. Introduction}\nl
A {\em Beurling generalized prime system}  ${\cal P}$ is an unbounded sequence of real numbers $p_1,p_2,p_3,\ldots$ satisfying
\[ 1<p_1\leq p_2\leq\cdots\leq p_n\leq \cdots . \]   
We call these numbers {\em generalized primes} (or $g$-primes), and from them we form the system ${\cal N}$ of {\em generalized integers} (or $g$-integers) associated to ${\cal P}$. These are the numbers of the form
\[ p_1^{\alpha_1}p_2^{\alpha_2}\cdots p_k^{\alpha_k}\]
where $\alpha_1, \ldots,\alpha_k\in\N_0$. 
In other words, ${\cal N}$ (viewed as a multi-set) is the semi-group generated by the (multi-set) ${\cal P}$ under multiplication. Such systems were first defined and investigated by Beurling \cite{B} in 1937 and have been studied by many researchers since then (see for instance  \cite{BD}, \cite{DZ} and the numerous references therein). Attached to these systems are the  counting functions 
\[ \pi_ \mathcal{P}(x) = \sum_{\tiny \begin{array}{c} p \le x \\ p\in {\cal{P}} \end{array}} 1,\quad N_ \mathcal{P}(x) = \sum_{\tiny \begin{array}{c} n \le x \\ n\in {\cal{N}} \end{array}} 1,\quad \psi_ \mathcal{P}(x) = \sum_{\tiny \begin{array}{c} p^k \leq x \\ p\in {\cal{P}}\\ k\in\N \end{array}} \log p,\]
which generalize the usual counting functions.
In each case, the sum is over all possible elements from the multi-set $\cal{P}$ or $\cal{N}$ with the given constraint. We are also interested in the generalized  {\em M\"{o}bius} function defined to be $\mu_{\cal{P}}(1)=1$, $\mu_ \mathcal{P}(p_{i_1}\cdots p_{i_k})=(-1)^k$ for distinct $g$-primes (i.e. $i_1,\ldots, i_k$ are distinct) and zero otherwise. Strictly speaking this need not be a function if two such products are numerically the same. In any case, we define the sum function 
\[M_\mathcal{P}(x) = \!\!  \sum_{\tiny \begin{array}{c} n \leq x \\ n\in {\cal{N}} \end{array}}  \mu_\mathcal{P}(n).\]
This generalizes the usual $M(x) = \sum_{n\le x} \mu(n)$. The associated {\em Beurling zeta function} is defined as usual by
\[ \zeta_ \mathcal{P}(s) = \prod_{p\in\mathcal{P}} \frac{1}{1-\frac{1}{p^s}} = \sum_{n\in\mathcal{N}} \frac{1}{n^s}.\]
We are interested in systems for which one or more of $\psi_ \mathcal{P}(x)-x$, $N_ \mathcal{P}(x)-\rho x$, or $M_ \mathcal{P}(x)$ is $O(x^\theta)$ for some $\theta<1$ (and $\rho>0$). More precisely, we define three numbers $\alpha,\beta,\gamma$ by the following:
\begin{align*} 
\psi_ \mathcal{P}(x) & = x+O(x^{\alpha+\eps})\tag{1}\\
N_ \mathcal{P}(x) & = \rho x+O(x^{\beta+\eps})\tag{2}\\
M_ \mathcal{P}(x) & = O(x^{\gamma+\eps})\tag{3}
\end{align*}
hold for all $\eps>0$ but no $\eps<0$. For example, for $\mathcal{N}=\N$, $\beta=0$ while $\alpha=\gamma\ge\frac{1}{2}$ due to the Riemann zeros. At the outset we are only interested in those systems for which the abscissa of convergence of the Dirichlet series for $\zeta_{\cal{P}}$ is 1. Thus $\alpha,\beta,\gamma\in [0,1]$ in any case.

For (1) and (2) to hold simultaneously for some $\alpha,\beta<1$ is akin to having a kind of Riemann Hypothesis being true for such a system. In \cite{Z2}, it was shown that such a system does exist with $\alpha,\beta\le\frac{1}{2}$. On the other hand, in \cite{H}, it was shown that it is impossible to have both $\alpha$ and $\beta$ less than $\frac{1}{2}$. 

\medskip

We note that (3) is related to an interesting problem in its own right: {\em how small can $M_ \mathcal{P}(x)$ be made for a system with abscissa\footnote{The abscissa of convergence of $\zeta_{\mathcal{P}}(s)$. With abscissa $\sigma_c$, we trivially have $M_{\mathcal{P}}(x)\ll x^{\sigma_c+\eps}$. Without the condition on the abscissa, $M_\mathcal{P}(x)$ can even be bounded: take $\mathcal{P} = \{ 2^{2^n}: n\in\N_0\}$. Then $M_\mathcal{P}(x)=0, 1$ or $-1$.}  equal to $1$?} In other words, how much cancellation can occur in the sum for $M_\mathcal{P}(x)$? Of course, for $\mathcal{N}=\N$, $M(x)=\Omega(\sqrt{x})$ on account of the Riemann zeros, but without this knowledge it is not clear how to even prove $M(x)=\Omega(x^a)$ for some $a>0$. This is similar to a question of Kahane and Saias \cite{KS} who ask how small $\sum_{n\le x} f(n)$ can be for $f$ completely multiplicative.

It is also related to the more general question of the size of $M_\mathcal{P}(x)$ and how it relates to the other functions. For example, much work has been done to determine under what conditions one has $M_\mathcal{P}(x)=o(x)$ (see for example Chapter 14 of \cite{DZ}). Zhang \cite{Z1} was the first to note that PNT is not equivalent to $M_\mathcal{P}(x)=o(x)$. For the most general results giving $M_\mathcal{P}(x)=o(x)$, see the very recent papers \cite{DDV} and \cite{DMV}. 

\medskip

Our main result is the following:\nl

\noindent
{\bf Theorem 1}\nl
{\em Of the numbers $\alpha,\beta,\gamma$, the two largest must be the same and at least $\frac{1}{2}$.}\nl

This result implies that for $M_ \mathcal{P}(x)=O(x^\gamma)$ with $\gamma<\frac{1}{2}$ to hold, we need the system to have somewhat chaotic $g$-primes and $g$-integers; i.e. the errors in (1) and (2) have to be $\Omega(x^{\frac{1}{2}-\eps})$ for every $\eps>0$. It may be conjectured that having $\gamma<\frac{1}{2}$ is actually impossible. 
\nl

\noindent
{\bf 2. Some Relevant Results}\nl
In order to prove the main result we shall need some relevant notions as well as existing results about $g$-prime systems. 

\medskip

Let $f(s) = \sum_{n=1}^\infty \frac{a_n}{b_n^s}$ be a generalized Dirichlet series where $b_n>0$  is strictly increasing with finite abscissae of absolute convergence $\sigma_a$. Suppose $f$ has a meromorphic continuation to $H_\alpha:=\{s\in\C: \Re s>\alpha\}$. We say $f$ has {\em finite order} in $H_\alpha$ if
\[ f(\sigma+it)\ll |t|^\lambda\qquad\mbox{($|t|\ge 1)$}\]
for $\sigma>\alpha$. As such, we can define the {\em Lindel\"{o}f function} $\mu_f(\sigma)$ to be the infimum of such $\lambda$. It is well-known that $\mu_f$ is non-negative, decreasing and convex and for $\sigma>\sigma_a$, $\mu_f(\sigma)=0$.\nl

The following result about such Dirichlet series and ``counting function''
\[ A(x) := \sum_{b_n\le x} a_n\]
 was essentially proved in \cite{H}, Proposition 3 (see also \cite{HL}, Theorem 2.1). It was proven for the case where $a_n\ge 0$ such that $a_n\ll n^\eps$ for all $\eps>0$. This latter condition however is not necessary. Also we shall require a particular case when $a_n$ is also sometimes negative. \nl

\noindent
{\bf Theorem A}\nl
{\em Let $f(s)=\sum_{n=1}^\infty \frac{a_n}{b_n^s}$ have abscissa of convergence $\sigma_c\le 1$. Suppose that for some $\alpha \in [0,1)$ and $c\in\C$, we have 
\[ A(x) = cx+ O(x^{\alpha+\eps}) \quad  \mbox{for all}\,\, \eps >0. \tag{4} \]
Then $f(s)$ has an analytic continuation to $H_\alpha\setminus\{1\}$ with a simple pole at $s=1$ with residue\footnote{Of course, if $c=0$, the pole is removable.} $c$ and $f$ has finite order; indeed $\mu_f(\sigma)\le 1$ for $\sigma>\alpha$.

Conversely, suppose that for some $\alpha\in [0,1)$, $f(s)$ has an analytic continuation to $H_\alpha$ except for a simple pole at $s=1$ with residue $c$. Further assume that 
$\mu_f(\sigma)=0$ for $\sigma>\alpha$ and either (i)  $a_n\ge 0$ or
\[ \hspace{-2in}\mbox{(ii)}\hspace{1.5in}\sum_{x-1<b_n\le x} |a_n| = O(x^{\alpha+\eps}) \quad  \mbox{for all}\,\, \eps >0. \tag{5} \]
Then $(4)$ holds.}\nl

\noindent
{\em Proof.}\, The proof of the first part is standard and follows on writing 
\[ f(s) = s\int_1^\infty \frac{A(x)}{x^{s+1}}\, dx = \frac{cs}{s-1} +s\int_1^\infty \frac{A(x)-cx}{x^{s+1}}\, dx,\]
and noting that the integral on the right converges absolutely to a holomorphic function on $H_\alpha$.

\medskip

For the converse, we follow the proof of Proposition 3 in \cite{H} as much as possible. This leads to
\[ A(x) - cx\ll \frac{x}{T^{1-\eps}}+x^{\alpha +\eps}T^{\eps}+\frac{x^{1+\eps}}{T}+\frac{x}{T}\sum_{\frac{x}{2}<b_n<2x}\frac{|a_n|}{|b_n- x|}\tag{6}\]
for every $T>1$ and $\eps>0$ --- see equation (3.7) of \cite{H}.

Now, as in \cite{H}, consider $x$ such that 
\[ \Bigl(x-\frac{1}{x^2}, x+\frac{1}{x^2}\Bigr)\cap \{ b_k:k\in\N\}=\emptyset.\tag{7}\]
For such $x$, $|b_n-x|\ge\frac{1}{x^2}$ for all $n$ and the sum on the right in (6) is at most
\[ x^2 \sum_{\frac{x}{2}<b_n<2x} |a_n|.\]
In case (i), this is $O(x^{3+\eps})$, while in case (ii), it is $O(x^{3+\alpha+\eps})$ by (5).

Taking $T=x^4$ in (6) shows that (4) holds whenever $x\to\infty$ satisfying (7). As shown in \cite{HL} (see (2.3)), for every $x$ there exist $x_1\in (x-1,x)$, $x_2\in (x,x+1)$ such that $x_1,x_2$ satisfy (7). Thus (4) holds for $x_1$ and $x_2$.

\medskip

For case (i), positivity of $a_n$ implies $A(x_1)\le A(x)\le A(x_2)$. Hence (4) follows for $x$.
  
For case (ii), we use (5). We have
\[ |A(x)-A(x_1)|\le 
\sum_{x-1<b_n\le x} |a_n| \ll x^{\alpha+\eps}\]
by (5). Hence (4) follows.
\bo

We also require the following result from \cite{HL} (Theorem 2.3).\nl

\noindent
{\bf Theorem B}\nl
{\em Suppose $(1)$ and $(2)$ hold for some $\alpha,\beta<1$. Then for $\sigma>\Theta:=\max\{\alpha,\beta\}$ and uniformly for $\sigma\ge\Theta+\delta$ (any $\delta>0)$,
\[ \frac{\zeta_{\mathcal{P}}^\prime(\sigma+it)}{\zeta_{\mathcal{P}}(\sigma+it)} = O\Bigl( (\log |t|)^{\frac{1-\sigma}{1-\Theta}+\eps}\Bigr)\quad\mbox{ and }\quad \zeta_{\mathcal{P}}(\sigma+it),\frac{1}{\zeta_{\mathcal{P}}(\sigma+it)} = O\Bigl(\exp\Bigl\{ (\log |t|)^{\frac{1-\sigma}{1-\Theta}+\eps}\Bigr\}\Bigr)\]
for all $\eps>0$. In particular, for $\sigma>\Theta$, the Lindel\"{o}f functions for $\zeta_{\mathcal{P}}$ and $\frac{1}{\zeta_{\mathcal{P}}}$ are zero.}\nl

Actually, the statement of Theorem 2.3 in \cite{HL} does not mention $\frac{1}{\zeta_{\mathcal{P}}}$ but the proof, which argues from $\log \zeta_{\mathcal{P}}$ clearly applies also to $-\log \zeta_{\mathcal{P}} = \log \frac{1}{\zeta_{\mathcal{P}}}$.\nl

Also, we have the following two consequences as described at the end of section 2 in \cite{HL}:
\begin{enumerate} 
\item If $\alpha>\beta$, then $\zeta_{\mathcal{P}}$ has infinitely many zeros on, or arbitrarily close to, the line $\sigma=\alpha$.
\item If $\alpha<\beta$, then $\zeta_{\mathcal{P}}$ and $\frac{1}{\zeta_{\mathcal{P}}}$ have infinite order in the strip $ \{ s \in \C: \alpha<\Re s<\beta \}$.
\end{enumerate}
\


\noindent
{\bf Proof of Theorem 1}\nl
Let $\Theta:=\max\{\alpha,\beta\}$. We use the converse part of Theorem A with $f(s)=\frac{1}{\zeta_{\mathcal{P}}(s)}$. This function has an analytic continuation to $H_\Theta$ and, by Theorem B, has zero order here. Further, $A(x)=M_{\mathcal{P}}(x)$ and 
\[  \sum_{\tiny \begin{array}{c}x-1< n \leq x \\ n\in {\cal{N}} \end{array}}  |\mu_\mathcal{P}(n)| \le N_{\mathcal{P}}(x)-N_{\mathcal{P}}(x-1)\ll x^{\beta+\eps}\le  x^{\Theta+\eps}.\] 
Thus (5), and hence (4), holds (with $c=0$). That is, $M_{\mathcal{P}}(x)=O(x^{\Theta+\eps})$; i.e. $\gamma\le\Theta$.

\medskip

Now suppose $\alpha>\beta$. Then $\zeta_{\mathcal{P}}$ has infinitely many zeros on, or arbitrarily close to, the line $\sigma=\alpha$. Thus $\gamma\ge\alpha-\delta$ for any $\delta>0$; i.e. $\gamma\ge\alpha$ and so $\gamma=\alpha$.

\medskip

Now suppose $\alpha<\beta$. Then the Lindel\"{o}f functions for $\zeta_{\mathcal{P}}$ and $1/{\zeta_{\mathcal{P}}}$ are infinite for $\sigma<\beta$. Thus we cannot have $\gamma<\beta$ by the first part of Theorem A with $A(x)=M_{\mathcal{P}}(x)$; i.e. $\gamma=\beta$. 

\medskip

Thus if $\alpha\ne\beta$, then $\gamma=\Theta$. Hence the two largest numbers are always equal. Finally, since $\max\{\alpha,\beta\}\ge\frac{1}{2}$, we see that in all three cases the largest pair is always at least $\frac{1}{2}$.\bo 

 \medskip

\noindent
{\bf 2. Systems with different $\alpha,\beta,\gamma$. }\nl
It is perhaps of interest to see if it really is possible that each of $\alpha,\beta$ or $\gamma$ can be strictly less than the other two and whether it can be less than $\frac{1}{2}$. 
\begin{enumerate}
\item $\beta<\alpha=\gamma$. For $\mathcal{N}=\N$, we have $\beta=0$ and, under the Riemann Hypothesis, $\alpha=\gamma=\frac{1}{2}$. Unconditionally, we only have $\alpha=\gamma=\Theta$ where $\Theta=\sup\{\Re \rho:\zeta(\rho)=0\}$.  
\item $\alpha<\beta=\gamma$. In the final discussion of \cite{H}, a $g$-prime system was given with $\alpha=0$. Namely, take $p_n=R^{-1}(n)$, where $R$ is the strictly increasing function on $[1,\infty)$ defined by
\[ R(x) = \sum_{k=1}^{\infty} \frac{(\log x)^k}{k!k\zeta(k+1)},\]
where $\zeta(\cdot)$ is the Riemann zeta-function. As such, one has $\psi_{\mathcal{P}}(x)= x+O(\log x\log\log x)$.
By Theorem 1, $\beta=\gamma$,  but what this common value is is not clear, except that it lies in $[\frac{1}{2},1]$. 
\item $\gamma<\alpha=\beta$. For this we can use the example $\mathcal{P}=\P\sqcup \P^{1/\beta}$ with $\beta\in (0,1)$. Using Dirichlet's hyperbola method, we have
\[N_{\mathcal{P}}(x) =\sum_{mn^{1/\beta}\le x}1 = \zeta\Bigl(\frac{1}{\beta}\Bigr)x + \zeta(\beta)x^{\beta}+O(x^{\frac{\beta}{1+\beta}})\]
(see \cite{H2} where this calculation was done). Furthermore, $\psi_{\mathcal{P}}(x) = \psi(x)+\psi(x^\beta) = x+x^\beta + O(x^{\frac{1}{2}+\eps})$ on RH. Thus $\alpha=\beta$. But,
with $M(x) = \sum_{n\le x} \mu(n)$,
\[M_{\mathcal{P}}(x)  =\sum_{mn^{1/\beta}\le x}\mu(m)\mu(n) = \sum_{n\leq a^{\beta}}M\Bigl(\frac{x}{n^{1/\beta}}\Bigr) + \sum_{n\leq b} M\Bigl( \Bigl(\frac{x}{n}\Bigr)^{\beta}\Bigr) - M(a^{\beta})M(b)\]
for any $ab=x$.  Putting $a=x^{\lambda}$ and using the bound $M(x)\ll x^{\frac{1}{2}+\eps}$ gives
\begin{align*}
M_{\mathcal{P}}(x) &\ll \sum_{n\leq x^{\lambda\beta}}\Bigl(\frac{x}{n^{1/\beta}}\Bigr)^{\frac{1}{2}+\eps} + \sum_{n\leq x^{1-\lambda}} \Bigl( \Bigl(\frac{x}{n}\Bigr)^{\beta}\Bigr)^{\frac{1}{2}+\eps} + (x^{\lambda\beta})^{\frac{1}{2}+\eps}x^{(1-\lambda)(\frac{1}{2}+\eps)}\\
& \ll \Bigl(x^{\frac{\beta}{2}+\lambda(1-\frac{\beta}{2})} + x^{\frac{1}{2}+(1-\lambda)(\beta-\frac{1}{2})} + x^{\frac{\lambda}{2}+(1-\lambda)\frac{\beta}{2}}\Bigr)x^\eps.
\end{align*}
Choosing $\lambda=\frac{\beta}{1+\beta}$ optimally shows that $M_{\mathcal{P}}(x)\ll x^{\frac{3\beta}{2(1+\beta)}+\eps}$ for all $\eps>0$. Thus $\gamma \le \frac{3\beta}{2(1+\beta)}<\beta$. Note that $\gamma\ge\frac{1}{2}$, since $\frac{1}{\zeta_{\mathcal{P}}(s)} = \frac{1}{\zeta(s)\zeta(s/\beta)}$ has poles on the $\frac{1}{2}$-line. 
\end{enumerate}

\medskip

\noindent
{\bf Open problems}\nl
1) From (a) and (b) above we have systems with $(\alpha,\beta,\gamma)=(a,0,a)$ and $(0,b,b)$ for some $a,b\in [\frac{1}{2},1]$. Can we find, unconditionally, such systems with $a<1$ and $b<1$?

\medskip

\noindent
2) In (c) above we have a system, conditional on RH, with $(\alpha,\beta,\gamma)=(c,c,d)$ with $\frac{1}{2}\le d<c<1$. Can we find one unconditionally, with $d<1$. Furthermore, can we find one with $d<\frac{1}{2}$?

\end{document}